\documentclass[12pt,a4paper]{amsart}
\usepackage{amsfonts}
\usepackage{amsthm}
\usepackage{amsmath}
\usepackage{amscd}
\usepackage[latin2]{inputenc}
\usepackage{t1enc}
\usepackage[mathscr]{eucal}
\usepackage{indentfirst}
\usepackage{graphicx}
\usepackage{graphics}
\usepackage{pict2e}
\usepackage{epstopdf}

\numberwithin{equation}{section}

\theoremstyle{plain}
\newtheorem{Th}{Theorem}[section]

 \theoremstyle{definition}
\newtheorem{Def}[Th]{Definition}

\newtheorem{Rem}[Th]{Remark}
\newtheorem{?}[Th]{Problem}

\newcommand{\ch}{\mathrm{ch}}

\begin{document}

\title{Co-adjoint polynomial}

\author{P\'eter Csikv\'ari} 
 \address{E\"{o}tv\"{o}s Lor\'{a}nd University \\ Department of Computer 
Science \\ H-1117 Budapest
\\ P\'{a}zm\'{a}ny P\'{e}ter s\'{e}t\'{a}ny 1/C \\ Hungary} 
\email{peter.csikvari@gmail.com}

 \thanks{The author  is partially supported by the National Science Foundation under grant no. DMS-1500219 and by the Hungarian National Research, Development and Innovation Office, NKFIH K109684, and by the ERC Consolidator Grant 648017.}


\begin{abstract} In this note we study a certain graph polynomial arising from
a special recursion. This recursion is a member of a family of four recursions where the other three recursions belong to the chromatic polynomial, the modified matching polynomial, and the adjoint polynomial. The four polynomials have many properties in common, for instance all of them are of exponential type, i. e., they satisfy the identity
$$\sum_{S\subseteq V(G)}f(G[S],x)f(G[V\setminus S],y)=f(G,x+y)$$
for every graph $G$.

It turns out that the new graph polynomial is a specialization of the Tutte polynomial.    
\end{abstract}

\maketitle

\section{Introduction} Throughout this paper all graphs are simple.  
Let us consider the following recursion for a graph
polynomial. Let $e=(u,v)\in E(G)$ and assume that $P(G,x)$ satisfies the
following recursion formula 
$$P(G,x)=P(G-e,x)-P(G\Delta e,x),$$
where $G\Delta e$ denotes the following graph. We delete the vertices $u$ and
$v$ from $G$ and replace it with a vertex $w$ which we connect to those
vertices of $V(G)-\{u,v\}$ that were adjacent to exactly one of $u$ and $v$
in $G$. In other words, we connect $w$ with the symmetric difference of
$N(u)\setminus\{v\}$ and  $N(v)\setminus\{u\}$. The $\Delta$ in the recursive formula refers to this
symmetric difference. Let $\overline{K_n}$ be the
empty graph on $n$ vertices and let $P(\overline{K_n},x)=x^n$. This completely
determines the graph polynomial $P(G,x)$ by induction on the number of
edges. On the other hand, it is not at all clear that this graph polynomial
exists since we can determine $P(G,x)$ by choosing edges in different orders
and we might not get the same polynomial. It will turn out that this polynomial
indeed exists and it is a specialization of the Tutte polynomial. Let us call
this graph polynomial co-adjoint polynomial for lack of a  better name.  
\medskip

What motivates this recursive formula of $P(G,x)$? Let us consider the
following three graph polynomials.  
\medskip

1. Let $M(G,x)=\sum_{k=0}^n(-1)^{k}m_k(G)x^{n-k}$ be the (modified) matching
polynomial \cite{god3,god4,hei} where $m_k(G)$ denotes the number of matchings
of size $k$ with the convention $m_0(G)=1$. Then $M(G,x)$ satisfies the
following recursive formula: 
let $e=(u,v)\in E(G)$ then
$$M(G,x)=M(G-e,x)-M(G\emptyset e,x)=M(G-e,x)-xM(G-\{u,v\},x),$$
where $G\emptyset e$ denotes the following graph. We delete the vertices $u,v$
from $G$ and replace it with a vertex $w$ which we do not connect with
anything. 
\medskip

2. Let $\ch(G,x)$ be the chromatic polynomial \cite{rea}. It is known that it
satisfies the following recursive formula. Let $e=(u,v)\in E(G)$ then 
$$\ch(G,x)=\ch(G-e,x)-\ch(G\cup e,x)=\ch(G-e,x)-\ch(G/e,x),$$
where $G/e=G\cup e$ denotes the following graph. We delete the vertices $u,v$
from $G$ and replace it with a vertex $w$ which we connect with the \emph{union} of
$N(u)\setminus\{v\}$ and  $N(v)\setminus\{u\}$.  
\medskip

3. Let $h(G,x)$ be the following graph polynomial. Let $a_k(G)$ be the number
of ways one can cover the vertex set of the graph $G$ with exactly $k$
disjoint cliques of $G$. Let 
$$h(G,x)=\sum_{k=1}^n(-1)^{n-k}a_k(G)x^k.$$
The graph polynomial $h(G,x)$ is called adjoint polynomial \cite{liu1,liu2}
(most often without alternating signs of the coefficients). Then $h(G,x)$
satisfies the following recursive formula. Let $e=(u,v)\in E(G)$ then 
$$h(G,x)=h(G-e,x)-h(G\cap e,x),$$ 
where $G\cap e$ denotes the following graph. We delete the vertices $u,v$ from
$G$ and replace it with a vertex $w$ which we connect with the \emph{intersection}
of $N(u)\setminus\{v\}$ and  $N(v)\setminus\{u\}$.  
\bigskip

\begin{figure}[h!]
\begin{center}
\scalebox{.65}{\includegraphics{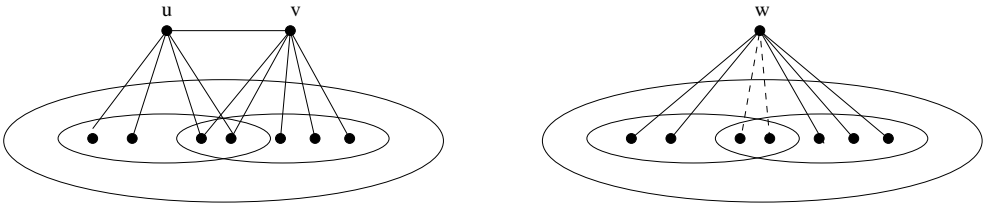}}
\caption{$f(G,x)=f(G-e,x)-f(G',x)$, where in $G'$ we consider the normal, dashed,
  all or no edges according to whether $f$ is the adjoint, co-adjoint, chromatic or
  matching polynomial, respectively.}   
\end{center}
\end{figure}

Now it is clear that the co-adjoint polynomial is the natural fourth member of
this family.  
\medskip

This paper is organized as follows. In the next section we prove that the
co-adjoint polynomial is a specialization of the Tutte polynomial, in particular it exits. The
third section is concerned with corollaries of this result. In the
last section we study the co-adjoint polynomials of complete graphs and
balanced complete bipartite graphs. 

\section{Specialization of the Tutte polynomial}

The Tutte polynomial of a graph $G=(V,E)$ is defined as follows.
$$T(G,x,y)=\sum_{A\subseteq E} (x-1)^{k(A)-k(E)}(y-1)^{k(A)+|A|-|V|},$$
where $k(A)$ denotes the number of connected components of the graph $(V,A)$.
 
In statistical physics one often studies the following form of the Tutte
polynomial:
$$Z_G(q,v)=\sum_{A\subseteq E}q^{k(A)}v^{|A|}.$$
The two forms are essentially equivalent:
$$T(G,x,y)=(x-1)^{-k(E)}(y-1)^{-|V|}Z_G((x-1)(y-1),y-1).$$
Both forms have several advantages. For instance, it is easy to generalize the
latter one to define the multivariate Tutte-polynomial. Let us assign a
variable $v_e$ to each edge $e$ and set
\begin{align} \label{multivariate Tutte}
Z_G(q,\underline{v})=\sum_{A\subseteq E}q^{k(A)}\prod_{e\in E}v_e.
\end{align}

Note that the chromatic polynomial of graph $G$ is
$$\ch(G,x)=Z_G(x,-1)=(-1)^{|V|-k(G)}x^{k(G)}T(G,1-x,0).$$
The main result of this section is the following.

\begin{Th} \label{Tutte} Let $G$ be a simple graph. Let $T(G,x,y)$  be the Tutte polynomial of $G$. Set
$$P(G,x)=\frac{1}{2^{|V(G)|}}Z_G(2x,-2)=(-1)^{|V|-k(G)}x^{k(G)}T(G,1-x,-1).$$
Then $P(G,x)$ satisfies the recursion formula
$$P(G,x)=P(G-e,x)-P(G\Delta e,x),$$
for an arbitrary edge $e$, and $P(\overline{K_n},x)=x^n$.
In particular, the polynomial defined by the above recursion indeed exists.
\end{Th}

\begin{Rem} The Tutte polynomial satisfies the following
  recursive formulas: 
$$T(G,x,y) = T(G-e,x,y)+T(G/e,x,y)$$
if e is neither a loop nor a bridge and 
$$T(G,x,y) = xT(G-e,x,y)$$
if $e$ is a bridge and
$$T(G,x,y) = yT(G/e,x,y)$$
if $e$ is a loop.

These formulas provide a straightforward way to prove Theorem~\ref{Tutte} by
induction. We will not follow this route since whenever we use these recursive
formulas we have to distinguish some cases according to the edge being a
bridge or not. After some steps the proof would split into too many cases.  
Instead we use the simpler form provided by the polynomial $Z_G(q,v)$.
\end{Rem}

\begin{proof} 
Let $|V(G)|=n$ and let us write
$$\frac{1}{2^n}Z(2x,-2)=\frac{1}{2^n}\sum_{A\subseteq
  E}(2x)^{k(A)}(-2)^{|A|}=\sum_{k=1}^n\left(\frac{1}{2^{n-k}}\sum_{A\subseteq E \atop
  k(A)=k}(-2)^{|A|}\right)x^k.$$ 
Set 
$$t_k(G)=\frac{1}{2^{n-k}}\sum_{A\subseteq E \atop k(A)=k}(-2)^{|A|}.$$
We need to prove that 
$$t_k(G)=t_k(G-e)-t_k(G\Delta e)$$
for an arbitrary edge $e$. So let us consider an arbitrary edge $e$. Clearly, in the
definition of $t_k(G)$, the sets $A$ not containing $e$ contribute a total of  $t_k(G-e)$ to the sum.

Now let us consider a set $A$ containing the edge $e$. Then one can consider $A-e$ as a set of edges in $G/e$ for which $k(A-e)=k$, whence
it follows that these sets contribute a total of $(-1)t_k(G/e)$ to the sum; note that
$|A-e|=|A|-1$, but $G/e$ has only $n-1$ vertices so the division and
multiplication by $2$ cancel each other, and only the term $-1$ remains from
the term $-2$. Hence
$$t_k(G)=t_k(G-e)-t_k(G/e).$$
Thus we only need to prove that
$$t_k(G\Delta e)=t_k(G/e).$$
So far we did not use anything about $G\Delta e$. Observe that $G\Delta e$ is
nothing else but the graph obtained from $G/e$ by deleting the multiple
edges. Let us consider a pair of multiple edges $e_1$ and $e_2$ incident with the same vertices. Assume that for
some edge set $A$ of $G/e$ not containing $e_1,e_2$ we have $k(A\cup
\{e_1\})=k$.  Then $k(A\cup \{e_2\})=k(A\cup \{e_1,e_2\})=k$ as well and they
contribute to the sum
$$(-2)^{|A\cup \{e_1\}|}+(-2)^{|A\cup \{e_2\}|}+(-2)^{|A\cup
  \{e_1,e_2\}|}=(-2)^{|A|}((-2)+(-2)+(-2)^2)=0.$$
Hence we can delete the multiple edges from $G/e$ without changing the value of
$t_k(.)$:
$$t_k(G/e)=t_k(G\Delta e).$$

Hence $t_k(G)=t_k(G-e)-t_k(G\Delta e)$ for an arbitrary edge $e$. Consequently,
$$P(G,x)=P(G-e,x)-P(G\Delta e,x).$$
To see that and $P(\overline{K_n},x)=x^n$ observe that for $k<n$ we have 
$$t_k(G)=\frac{1}{2^{n-k}}\sum_{A\subseteq E \atop k(A)=k}(-2)^{|A|}=0$$
and for $k=n$ we have $t_n(G)=1$.
This completes the proof.

\end{proof}
 
\begin{Rem} By the recursive formula
$$P(G,x)=P(G-e,x)-P(G\Delta e,x)$$
it is easy to prove that the coefficients $t_k(G)$ have alternating signs.
On the other hand, it is not clear from the expressions given for $t_k(G)$.  
\end{Rem}

\begin{Rem} A surprising corollary of Theorem~\ref{Tutte} is that $|P(G,1)|=0$
  or $1$ and it is $1$ if and only if the graph is Eulerian, i.e., all degrees
  are even. This follows from the fact that $|T(G,0,-1)|$ counts the nowhere-$0$
  $\mathbb{Z}_2$-flows (note that the flow polynomial is also a specialization
  of the Tutte polynomial), and a  nowhere-$0$ $\mathbb{Z}_2$-flow is simply a
  flow taking the value $1$ on all edges, this immediately implies the claim.
\end{Rem}

\section{Exponential type graph polynomials}

In the introduction we considered four graph polynomials: the matching
polynomial, the chromatic polynomial, the adjoint polynomial and our new graph
polynomial, the co-adjoint polynomial. Surprisingly, they all belong to a very
special class of graph polynomials, the so-called exponential type graph
polynomials.  

\begin{Def} We say that the graph polynomial $f$ is of exponential type if
for every graph $G=(V(G),E(G))$ we have $f(\emptyset,x)=1$  and $f(G,x)$
satisfies that
$$\sum_{S\subseteq V(G)}f(G[S],x)f(G[V\setminus S],y)=f(G,x+y),$$
where $G[S]$ denotes the induced subgraph of $G$ on vertex set $S$. 
\end{Def}

Note that Gus Wiseman \cite{wis} calls these graph polynomials
\textit{binomial-type}. 

One can deduce from the definition that the chromatic polynomial is of
exponential type. For the matching polynomial and the adjoint polynomial this
follows from Theorem~\ref{exp-b} below. This is a structure theorem for the
exponential type graph polynomials proven in \cite{fre}. For the co-adjoint polynomial this is
simply a special case of the following much more general statement.

\begin{Th} \cite{sco-sok2} For the multivariate Tutte-polynomial
  $Z_G(q,\underline{v})$ we have 
$$\sum_{S\subseteq V(G)}Z_{G[S]}(q_1,\underline{v})Z_{G[V\setminus
      S]}(q_2,\underline{v})=Z_G(q_1+q_2,\underline{v}).$$ 
\end{Th}

The following theorem characterizes exponential type graph polynomials, see Theorem 5.1 of \cite{fre}.

\begin{Th} \cite{fre} \label{exp-b} Let $b$ be a function from the class of graphs to the
  complex numbers. Let us define the graph polynomial $f_b$ as follows. Set
$$a_k(G)=\sum_{\{S_1,S_2,\dots ,S_k\}\in \mathcal{P}_k}b(S_1)b(S_2)\dots
b(S_k),$$
where the summation goes over the set $\mathcal{P}_k$ of all  partitions of
$V(G)$ into exactly $k$ non-empty sets. Then let 
$$f_b(G,x)=\sum_{k=1}^na_k(G)x^k,$$
where $n=|V(G)|$. Then
\medskip

(a) For any function $b$, the graph polynomial $f_b(G,x)$ is of exponential type.
\medskip

(b) For any exponential type graph polynomial $f$, there exist a graph function
$b$ such that $f(G,x)=f_b(G,x)$. More precisely, if $b(G)$ is the coefficient
of $x^1$ in $f(G,x)$ then $f=f_b$. 
\end{Th}

\begin{Rem} For the matching polynomial take $b_m(K_1)=1$ $b_m(K_2)=-1$ and
  $b_m(H)=0$ otherwise. For the adjoint polynomial consider
  $b_h(K_n)=(-1)^{n-1}$ for complete graphs $K_n$ and $b_h(H)=0$
  otherwise. This proves that the matching and the adjoint polynomials are
  indeed of exponential type.  
\end{Rem}

\begin{Rem} By the method of Alan Sokal \cite{sok} one can prove that the root
  of $P(G,x)$ of largest modulus has absolute value at most $KD$ where $D$ is
  the maximum degree of $G$, and  
$$K:=\inf_a \frac{a+e^a}{\log(1+ae^{-a})}\approx 7.963907.$$
Alan Sokal \cite{sok} proved this statement for the chromatic polynomial, and in fact, he proved this result for multivariate Tutte-polynomials where $v_e$ satifies $|1+v_e|\leq 1$ for each edge $e$, see Corollary 5.5 in the paper \cite{sok}. 

 This more general result covers the co-adjoint polynomial as we have $v_e=-2$. Note that this result would give an upper bound $2KD$ for the first sight, but since $P(G,x)=2^{-|V(G)|}Z(2x,-2)$ we can immediately get back the factor $2$ thereby providing an upper bound $KD$.
 Alternatively, Theorem 1.6 of \cite{fre} or results from the paper \cite{JPS} provide a weaker, but still linear bounds.
\end{Rem}
 
\section{Complete graphs and balanced complete bipartite graphs}

In this section we give the co-adjoint polynomial of some small graphs.
$$P(K_1,x)=x$$
$$P(K_2,x)=x^2-x$$
$$P(K_3,x)=x^3-3x^2+x$$
$$P(K_4,x)=x^4-6x^3+7x^2-2x$$
$$P(K_5,x)=x^5-10x^4+25x^3-20x^2+5x$$
$$P(K_6,x)=x^6-15x^5+65x^4-105x^3+70x^2-16x$$
$$P(K_7,x)=x^7-21x^6+140x^5-385x^4+490x^3-287x^2+61x$$
$$P(K_8,x)=x^8-28x^7+266x^6-1120x^5+2345x^4-2548x^3+1356x^2-272x$$
\bigskip

Clearly, the coefficient of  $x^1$ in $P(G,x)$ is
$(-1)^{|V|-1}T(K_n,1,-1)$ by Theorem~\ref{Tutte}. It is known that
$a_n=T(K_n,1,-1)$ counts the number of alternating permutations on $n-1$
elements. Let 
$(-1)^nP(K_n,-x)=p_n(x)$. The graph polynomial $P(G,x)$ is of exponential type.
Applying this observation to the complete graphs we obtain that
$$\sum_{k=0}^n{n \choose k}p_k(x)p_{n-k}(y)=p_n(x+y).$$
Hence the polynomials $(p_k(x))$ are of binomial type and consequently, we know that
$$\sum_{n=0}^{\infty}p_n(x)\frac{z^n}{n!}=\exp(xF(z)),$$
where 
$$F(z)=\sum_{n=1}^{\infty}a_n\frac{z^n}{n!}.$$
The exponential generating functions of the alternating permutations is known,
we only need to integrate it since the coefficients are translated:
$$F(z)=\int \frac{1+\sin z}{\cos z}=\ln \frac{1+\sin z}{\cos^2 z}.$$
\bigskip

For balanced complete bipartite graphs we have
$$P(K_{1,1},x)=x^2-x$$
$$P(K_{2,2},x)=x^4-4x^3+6x^2-2x$$
$$P(K_{3,3},x)=x^6-9x^5+36x^4-66x^3+51x^2-13x$$
$$P(K_{4,4},x)=x^8-16x^7+120x^6-488x^5+1112x^4-1360x^3+808x^2-176x$$
$$P(K_{5,5},x)=x^{10}-25x^9+300x^8-2100x^7+9150x^6-25030x^5+$$
$$+42020x^4-41020x^3+20785x^2-4081x$$ 
\medskip

The sequence of the coefficients of $x^1$ seems to be very interesting. Note
that not  only these numbers are $1,2,13,176,4081,\dots$, but the values of
$P(K_{n,n},-1)$ are also  these numbers. The same phenomenon occurs at \\
$(-1)^nP(K_n,-1)$ and the coefficients of $P(K_{n+2},x)$. In fact, these are
known results. The latter is a result of Merino, and both of these facts are 
special cases of the main result of \cite{GMMN} which asserts that under
certain conditions we have
$$T(G,1,-1)=T(G-\{u,v\},2,-1).$$
We do not give the conditions of their theorem here, but we note that
the complete graphs and complete bipartite graphs satisfy the conditions if
$(u,v)$ is an edge. 
\bigskip

\textbf{Acknowledgment.} We are very grateful to Mikl\'os B\'ona for various useful comments. We are also very grateful to the authors of \cite{GMMN} for including a table about $T(K_{m,n},2,-1)$ into their paper, which was crucial for us to make the right guess  about the studied graph polynomial.

Finally, a very big thanks go to the referees of this paper. I always appreciate the work of referees as it is a clearly volunteering work, but this time I am even more grateful to them for their numerous observations and advices.

\bibliographystyle{plain}

\bibliography{co-adjoint_bib}

\end{document}